# Possible Worlds, Incompleteness and Undefinability

This short squib looks at how using a broader definition of Gödel numbering to mimic the accessibility relation between possible worlds results in two-world systems that sidestep undecidable sentences as well as the Liar paradox.

## 1. Introduction

Possible worlds can be defined to be anything from objects, to variables to sets of sentences. As sets of sentences, they are maximally consistent in that, for every sentence, either it or its negation is in the set, but not both (maximality), and no finite subset of sentences can be used to derive a contradiction (consistency) (Hughes and Cresswell 2004:37-38). We can extend this definition to formulas in general, such that for any formula F, either F or ¬F is in the set, but not both.

There is also an accessibility relation R between worlds that indicates which worlds each world has access to. If a world has access to another one, then it has access to the sentences in that world, which enables it to have propositions concerning such sentences. This is precisely how the truth values of possible propositions (◊p) and necessary propositions (□p) at a given world $w_0$ are defined in modal logic:

$V(◊p, w_0) = 1 \leftrightarrow \exists w\ (w_0 R w \land V(p, w) = 1)$

$V(□p, w_0) = 1 \leftrightarrow \forall w\ (w_0 R w \rightarrow V(p, w) = 1)$

where V(p, w) is the truth valuation function of p at world w, 1 and 0 are the truth values corresponding to True and False, respectively, and $w_0 R w$ means that $w_0$ has access to w. What the above two definitions express is that if a sentence is true in at least one world that can be accessed by $w_0$, then it is possibly true at $w_0$, and if it is true in all worlds that can be accessed by $w_0$, then it is necessarily true at $w_0$ (Hughes and Cresswell 2004:38).

A second example, cases of which are the only ones that will be considered in this squib, consists of a world $w_k$ and one of its sentences $\varphi_k$, short for "*$\varphi$-in-$w_k$*", that happens to be true in $w_k$; symbolically, we could write this as $T_k(\mathbf{g_k(\varphi_k)})$ where $T_k$ is the truth predicate for $w_k$ and $\mathbf{g_k(\varphi_k)}$ is the numeral of the Gödel number $g_k(\varphi_k)$ where $g_k$ is the Gödel numbering used in association with $w_k$. Moreover, this sentence is also true in all worlds that have access to $w_k$ (Plantinga 1982:55), that is to say, $\forall i\ (w_i\ R\ w_k \rightarrow T_i(\mathbf{g_i(\varphi_k)}))$.

An illustration of this would be two worlds $w_1$ and $w_2$, such that the sky is blue in $w_1$, while it is plaid in $w_2$. Let $w_1\ R\ w_2$ and let $\varphi_2$ be the sentence "*the sky is plaid in $w_2$*". This would be true in $w_2$, so we have $T_2(\mathbf{g_2(\varphi_2)})$: "*it is true in $w_2$ that the sky is plaid in $w_2$.*". However, since $w_1\ R\ w_2$, we also have $T_1(\mathbf{g_1(\varphi_2)})$: "*it is true in $w_1$ that the sky is plaid in $w_2$.*"

There are probably many ways to indicate that a given expression is from a given world, but for the purposes of illustration, and for the sake of simplicity, one direct method is to append the number of a world to all expressions in that world. Thus, $(0 = 0)_2$ is an expression from $w_2$.

As number-theoretical statements like $0 = 0$ hold the same in all worlds, how does one world differ from another? It is due to predicates, like truth and provability predicates, whose domains consist of sets of numerals corresponding to the Gödel numbers of sentences.[1] These predicates vary from world to world as their domains depend on the accessibility relation. For example, if $w_1$ only has access to itself and to $w_2$, then the domain of its provability predicate $P_1$ is based on the sentences of $w_1 \cup w_2$. If $w_2$ only has access to $w_3$, then the domain of its predicate $P_2$ is just based on $w_3$. As a result, $P_1$ is distinct from $P_2$, and $w_1$, containing sentences that also contain $P_1$, is distinct from $w_2$.

## 2. Gödel numbering and worlds

Just as there are many ways to distinguish expressions from different worlds, there are many possible Gödel numberings to encode them. For the purposes of illustration, a slight adaptation of the Gödel numbering from Boolos, Burgess and Jeffrey[2], is used: "the scheme for coding finite sequences of numbers by single numbers based on prime decomposition" (2007:193):

| Symbol | ( | ) | , | ¬ | ∨ | ∃ | = | $v_i$ | $A_i^n$ | $f_i^n$ |
|---|---|---|---|---|---|---|---|---|---|---|
| Code | 1 | 3 | 5 | 7 | 9 | 11 | 13 | $2 \bullet 5^i$ | $2^2 \bullet 3^n \bullet 5^i$ | $2^3 \bullet 3^n \bullet 5^i$ |

Expressions, which are just strings of symbols, are encoded by taking the primes in increasing order and assigning the number code of each symbol as an exponent of the corresponding prime in the sequence of primes; the exceptions are the first and second primes in the sequence. The exponent of the first prime is the number designating the world to which the expression belongs; the exponent of the second prime is the length of the string that makes up the expression. Finally, all of the exponentialized primes are multiplied together to give the Gödel number of the string. As a result, the Gödel number for a string $x_0 \ldots x_{n-1}$ of length n in $w_k$ has the form $2^k 3^n 5^{g_k(x_0)} \ldots \pi(n+1)^{g_k(x_{n-1})}$ where the exponent of two is the value k from $w_k$ (the world number k), the exponent of three is the length of the string, of value n, and the exponents of the following primes are the Gödel numbers of the symbols of the string in order of appearance: $g_k(x_0)$ to $g_k(x_{n-1})$. Finally, as two is the zeroth prime and three is the first, coding the string requires going from five, the second prime, up to the $(n+1)^{st}$ prime, which is designated as $\pi(n+1)$.

We will now look at $(0 = 0)_2$ as an example [3]. Following Boolos et al. (2007:193), we start with the plain expression $0 = 0$, which is encoded as follows. First, the proper way of representing the expression is with the predicate first, followed by its arguments in parentheses, separated by commas, i.e., as =(0,0). Next, the number code for zero is that corresponding to $f_0^0$, which is $2^3 \bullet 3^0 \bullet 5^0 = 8$. So, the sequence of numbers corresponding to the symbols in the string are 13, 1, 8, 5, 8 and 3. The world number is 2 and the length of the string is 6, so the Gödel number is then $2^2 \bullet 3^6 \bullet 5^{13} \bullet 7^1 \bullet 11^8 \bullet 13^5 \bullet 17^8 \bullet 19^3$.

---

1 From here on, predicates whose domains are the numerals of Gödel numbers of sentences will be said to act on sentences or to have sentence arguments.
2 There are many other numberings available, such as that used by Smullyan (1992).
3 Of course, $(0 = 0)_2$ should really be written as $=_2(0,0)$, following the format $P_k(x,y)$ for a predicate from $w_k$.

As a result, it is easy to determine the world to which an expression $E$, with Gödel number $e$, belongs as well as its length ( $lh(e)$ ):

world(e) = lo(e,2), which is the exponent of the prime factor 2 of e, and
lh(e) = lo(e,3), which is the exponent of the prime factor 3 of e.

where lo(x,y) is the logarithm function defined in Boolos et al. (2007:79) and can be given as

$lo(x,y) = z \leftrightarrow (x>1) \wedge (y>1) \wedge ( (x \bmod y^z = 0) \rightarrow \forall t<x ( (x \bmod y^t = 0) \rightarrow t \leq z ) )$ ;
otherwise, $lo(x,y) = 0$.

## 3. Gödel numbering and the accessibility relation

Now, if a world $w_i$ has access to $w_k$, then the predicates of $w_i$ that take sentence arguments should be able to act on sentences from $w_k$. To this end, for each world $w_i$, let its corresponding Gödel numbering $g_i$ be defined as

$$(1) \quad g_i = \begin{cases} g_i^0(w_0) \\ \vdots \\ g_i^n(w_n) \\ \vdots \end{cases}$$

where the $g_i^k(w_k)$ are component Gödel numberings of $g_i$, which may be distinct from each other, so that for each world $w_k$, there is a corresponding component function of $g_i$ that acts on it. Furthermore, $0 \leq k \leq |W|$ where W is the set of all worlds. Eq. (1) can also be expressed as

$g_i(w_k) = g_i^k(w_k)$ over all k.

This generalization of the Gödel numbering function provides a way to express the accessibility relation. Suppose $w_i$ has no access to $w_k$, i.e., $\neg w_i R w_k$, then, as $w_i$ cannot access the formulas of $w_k$, they cannot be Gödel-numbered using $g_i$, i.e., $g_i^k(w_k) = \varnothing$. Conversely, if $g_i^k(w_k) = \varnothing$, then the formulas of $w_k$ cannot be arguments of any of the predicates of $w_i$, which would also be the case if $\neg w_i R w_k$. Therefore, we can make the following definition: $\neg w_i R w_k \leftrightarrow g_i^k(w_k) = \varnothing$. As a result, accessibility can be related to Gödel numbering in the following manner

**(2)**      $w_i \, R \, w_k \leftrightarrow g_i^k(w_k) \neq \varnothing$

Of special interest is the case where $w_i$ cannot access itself: if there is no encoding of its own formulas, it cannot have an undecidable Gödel sentence or a Liar paradox. These can only occur when R is reflexive since expressions like $\neg P_k(g_k(\varphi_k)) \leftrightarrow \varphi_k$ would only then be possible. Otherwise, and assuming at least $w_i \, R \, w_k$, we would have hybrid expressions of the form $\neg P_i(g_i(\varphi_k)) \leftrightarrow \varphi_k$, where $\neg P_i(g_i(\varphi_k))$ is a $w_i$-sentence in $w_i$, while $\varphi_k$ is a $w_k$-sentence in $w_k$, so that the whole expression cannot be a sentence of either world. This will also avoid sentences that claim not to be provable in either $w_i$ or $w_k$ as such sentences would be hybrid expressions.

**Regarding truth predicates:** Let (1) $w_i \, R \, w_k$, (2) $T_i(g_i(S_i))$ not hold for any sentence $S_i$ in $w_i$, the result of $\neg \, w_i \, R \, w_i$ or some other condition, and (3) $\Sigma_k$ hold, or be true, in $w_k$, then

**(3)**      $T_i(g_i(\Sigma_k))$

The advantage of this is that Tarski's Convention T does not hold as $T_i(g_i(\Sigma_k)) \leftrightarrow \Sigma_k$ is a hybrid expression. Thus, the Liar paradox is avoided, and there can be a truth predicate in $w_i$ for sentences in $w_k$. If we further require the same for $w_k$, then this two-world system is like a language-metalanguage system but without the infinite hierarchies of metalanguages.

**Regarding undecidable Gödel sentences:** If $w_i \, R \, w_k$, then an undecidable Gödel sentence $G_k$ is provable in $w_i$ by stipulating that it satisfy $Axiom_i(S)$. If (1) $w_i \, R \, w_k$ and (2) $T_i(g_i(S_i))$ does not hold for any $S_i$, then $T_i(g_i(G_k))$ since the fact that $G_k$ is always true in $w_k$ satisfies condition 3 of Eq. (3).

As an aside, an alternative to possible worlds is the use of copies of a theory of arithmetic. In this case, the accessibility relations are entirely replaced by their expression in terms of Gödel numberings, and one copy of the theory differs from another in that the domain of any given predicate that acts on sentences differs from copy to copy as a result of the distinct Gödel numberings associated with each copy.

## 4. Various two-world systems

The table in this section lists all of the ten two-world system types, consisting of $w_i$ and $w_k$, with all of the possible types of accessibility relation between them, the corresponding conditions on g and the corresponding possible predicates with sentence arguments. For each type, the table also indicates whether there are undecidable Gödel sentences (symbolized as $G_i$ or $G_k$), Liar paradoxes (symbolized as $L_i$ or $L_k$) or truth predicates (symbolized as $T_i$ or $T_k$). Only ten types are given instead of the full sixteen cases as the missing six cases are easily obtained by interchanging sub-indices i and k, and the results do not introduce any new forms, which explains why only $w_i \, R \, w_k$ appears and not $w_k \, R \, w_i$. Note that in cases 1 to 3, R is not reflexive, so all self-reference is avoided, resulting in no undecidable sentences and no Liar paradoxes. However, cases 2 and 3 allow for paradox-free truth predicates as in Eq. (3). There are also no undecidable Gödel sentences as axioms in cases 4-6, and $G_i$ is not an axiom of $w_k$ in case 9, as the necessary access between worlds is not present in these cases.

Some truth predicates may need to be suppressed entirely (e.g., in case of $L_i$, $T_i$ must be eliminated), or partially, e.g., in case of both $L_i$ and $T_i(S_k)$ for all $S_k$, the suppression of all $S_i$ from the domain of $T_i$ is required, leaving only all the $S_k$-type sentences, in order to avoid any $\neg T_i(S_i)$.

|   | **Types of $R^4$** | **Corresponding non-empty $g_i^k$** | **Predicates with S. A. [5] for all $S_i$, $S_k$** | **Undecidable Gödel sentences (G), Liars (L) and Truth predicates (T) acting on all $S_i$, $S_k$, $\Sigma_i$ and $\Sigma_k$ - type sentences** |
|---|---|---|---|---|
| 1 | None | None | None | None |
| 2 | $(w_i\ R\ w_k)$ | $g_i^k(w_k) \neq \varnothing$ | $P_i(\mathbf{S_k})$ [6] | $T_i(\mathbf{\Sigma_k})$<br>& No $G_i$, No $G_k$, No Liars |
| 3 | $(w_i\ R\ w_k)$<br>$(w_k\ R\ w_i)$ | $g_i^k(w_k) \neq \varnothing$<br>$g_k^i(w_i) \neq \varnothing$ | $P_i(\mathbf{S_k})$<br>$P_k(\mathbf{S_i})$ | $T_i(\mathbf{\Sigma_k})$<br>$T_k(\mathbf{\Sigma_i})$<br>No $G_i$, No $G_k$, No Liars |
| 4 | $(w_i\ R\ w_i)$ | $g_i^i(w_i) \neq \varnothing$ | $P_i(\mathbf{S_i})$ | $G_i$ unprovable & ($L_i \Rightarrow$ eliminate $T_i$) |
| 5 | $(w_i\ R\ w_i)$<br>$(w_k\ R\ w_k)$ | $g_i^i(w_i) \neq \varnothing$<br>$g_k^k(w_k) \neq \varnothing$ | $P_i(\mathbf{S_i})$<br>$P_k(\mathbf{S_k})$ | $G_i$ & $G_k$ both unprovable<br>$L_i$ & $L_k \Rightarrow$ eliminate both $T_i$ & $T_k$ |
| 6 | $(w_i\ R\ w_i)$<br>$(w_i\ R\ w_k)$ | $g_i^i(w_i) \neq \varnothing$<br>$g_i^k(w_k) \neq \varnothing$ | $P_i(\mathbf{S_i})$<br>$P_i(\mathbf{S_k})$ | $G_i$ unprovable & ($L_i \Rightarrow$ eliminate $T_i(S_i)$ )<br>$T_i(\mathbf{\Sigma_k})$, and due to $L_i \Rightarrow T_i$ is partial. |
| 7 | $(w_k\ R\ w_k)$<br>$(w_i\ R\ w_k)$ | $g_k^k(w_k) \neq \varnothing$<br>$g_i^k(w_k) \neq \varnothing$ | $P_k(\mathbf{S_k})$<br>$P_i(\mathbf{S_k})$ | $G_k$ & ($L_k \Rightarrow$ eliminate $T_k$)<br>Axiom$_i(\mathbf{G_k})$ ∴ $G_k$ is provable in $w_i$.<br>$T_i(\mathbf{\Sigma_k})$ ∴ $T_i(\mathbf{G_k})$; $G_k$ is true in $w_i$. |
| 8 | $(w_i\ R\ w_i)$<br>$(w_i\ R\ w_k)$<br>$(w_k\ R\ w_i)$ | $g_i^i(w_i) \neq \varnothing$<br>$g_i^k(w_k) \neq \varnothing$<br>$g_k^i(w_i) \neq \varnothing$ | $P_i(\mathbf{S_i})$<br>$P_i(\mathbf{S_k})$<br>$P_k(\mathbf{S_i})$ | $G_i$ & ($L_i \Rightarrow$ eliminate $T_i(S_i)$ )<br>$T_i(\mathbf{\Sigma_k})$, and due to $L_i \Rightarrow T_i$ is partial.<br>Axiom$_k(\mathbf{G_i})$ ∴ $G_i$ is provable in $w_k$.<br>$T_k(\mathbf{\Sigma_i})$ ∴ $T_k(\mathbf{G_i})$; $G_i$ is true in $w_k$. |
| 9 | $(w_i\ R\ w_i)$<br>$(w_k\ R\ w_k)$<br>$(w_i\ R\ w_k)$ | $g_i^i(w_i) \neq \varnothing$<br>$g_k^k(w_k) \neq \varnothing$<br>$g_i^k(w_k) \neq \varnothing$ | $P_i(\mathbf{S_i})$<br>$P_k(\mathbf{S_k})$<br>$P_i(\mathbf{S_k})$ | $G_i$ unprovable & ($L_i \Rightarrow$ eliminate $T_i(S_i)$ )<br>$G_k$ & ($L_k \Rightarrow$ eliminate $T_k$)<br>$T_i(\mathbf{\Sigma_k})$, and due to $L_i \Rightarrow T_i$ is partial.<br>Axiom$_i(\mathbf{G_k})$ ∴ $G_k$ is provable in $w_i$.<br>$T_i(\mathbf{\Sigma_k})$ ∴ $T_i(\mathbf{G_k})$; $G_k$ is true in $w_i$. |
| 10 | $(w_i\ R\ w_i)$<br>$(w_k\ R\ w_k)$<br>$(w_i\ R\ w_k)$<br><br>$(w_k\ R\ w_i)$ | $g_i^i(w_i) \neq \varnothing$<br>$g_k^k(w_k) \neq \varnothing$<br>$g_i^k(w_k) \neq \varnothing$<br><br>$g_k^i(w_i) \neq \varnothing$ | $P_i(\mathbf{S_i})$<br>$P_k(\mathbf{S_k})$<br>$P_i(\mathbf{S_k})$<br><br>$P_k(\mathbf{S_i})$ | $G_i$ & ($L_i \Rightarrow$ eliminate $T_i(S_i)$) $\Rightarrow T_i$ is partial.<br>$G_k$ & ($L_k \Rightarrow$ eliminate $T_k(S_k)$) $\Rightarrow T_k$ is partial.<br>$T_i(\mathbf{\Sigma_k})$ ∴ $T_i(\mathbf{G_k})$; $G_k$ is true in $w_i$.<br>Axiom$_i(\mathbf{G_k})$ ∴ $G_k$ is provable in $w_i$.<br>$T_k(\mathbf{\Sigma_i})$ ∴ $T_k(\mathbf{G_i})$; $G_i$ is true in $w_k$.<br>Axiom$_k(\mathbf{G_i})$ ∴ $G_i$ is provable in $w_k$. |

---

[4] In case 3, R is symmetric, in case 5 it is reflexive and in case 10 it is symmetric, reflexive and thus also transitive.
[5] S. A. stands for "Sentence Arguments".
[6] For simplicity, symbols in bold represent the numeral of the Gödel number of the symbol, so that $\mathbf{S_k} = \overline{g_k(S_k)}$.

## 5. Observations

1. Cases 1-3 are of interest as they avoid self-reference altogether, thus avoiding both the Liar and undecidable sentences, and cases 2 and 3 have truth predicates that are complete. However, the lack of self-reference is a drawback for these cases.

2. Cases 4 and 5 are not of interest as unprovable undecidable sentences and Liars hold in both.

3. Case 6 has an unprovable undecidable sentence and a partial truth predicate $T_i$, for which all $S_i$ must be suppressed from its domain to avoid the Liar paradox. Such partial predicates are not desirable as they are the result of an ad hoc suppression of a subset of their domain.

4. In case 7, $w_k$ has an undecidable Gödel sentence $G_k$ that is provable in $w_i$ since it can be stipulated to satisfy $Axiom_i(\mathbf{S})$, given that $w_i\ R\ w_k$, which allows for the encoding of $w_k$-sentences by $g_i$. $G_k$ is also true in $w_i$ as (i) $w_i\ R\ w_k$, (ii) $\neg w_i\ R\ w_i$, with the result that $T_i(\mathbf{g_i(S_i)})$ cannot hold for any $S_i$, and (iii) the fact that $G_k$ is always true in $w_k$ (undecidable Gödel sentences are always true in the system they reference). These three conditions satisfy those required by Eq. (3), the definition of the truth predicate. An advantage of case 7 is that $w_i$ cannot have its own undecidable sentence since it cannot access itself, and its truth predicate can only act on sentences from $w_k$. $T_k$, the truth predicate in $w_k$, however, is subject to the Liar paradox and has to be suppressed as in Tarski's hierarchy of languages, yet nothing is lost as $T_i$ predicates the truth of sentences in $w_k$, so $T_k$ is irrelevant to the system. The world $w_i$, thus acts as a kind of appendix to $w_k$, so that the two-world system of $w_i$ and $w_k$ contains all the true and provable sentences of $w_k$.

5. Case 8 is just like case 7 but with a true and provable $G_i$ in $w_k$; furthermore, unlike case 7, there is no need to eliminate any truth predicates as it has both a complete truth predicate $T_k$, and also a partial one $T_i$, which in any case is irrelevant to the system as only $T_k$ is needed for the truth of $G_i$.

6. Case 9 is particularly problematic as it has an unprovable undecidable Gödel sentence $G_i$, a Liar paradox in $w_k$, leading to the elimination of $T_k$, and a partial truth predicate $T_i$.

7. Finally, case 10 has two undecidable Gödel sentences, but both are true and provable in their respective opposite worlds; however, there are two partial truth predicates, one for each world.

8. Stacking of predicates exists in various forms in all cases save the first two.

Of the ten types of cases, the most interesting are 7, 8 and 10. The shortcomings of each case are as follows. Case 7 requires the elimination of a truth predicate to avoid the Liar; however, this predicate is irrelevant to the system. Case 8 has a partial truth predicate; however, this predicate is also irrelevant to the system. Which is worse depends on whether eradicating a truth predicate is worse than allowing for a partial one. Finally, case 10 has two partial truth predicates, both required by the system, which is a disadvantage due to their ad hoc nature.

## 6. Final note

It results that, by using a broader type of Gödel numbering function, the accessibility relation between worlds can be mimicked, and cases can be found that sidestep the limiting theorems of both Gödel and Tarski. This applies to systems that have theories of arithmetic contained in possible worlds, where worlds are sets of sentences, or to a set of copies of a theory of arithmetic where copies differ due to different versions of predicates that have sentence arguments since the domains of the different versions vary according to the different Gödel numberings assigned to each copy.